\documentclass[12pt,a4paper,twoside]{article}

\newcommand{\pageformat}[6]{\setlength{\hoffset}{-1in}
                  \setlength{\voffset}{-1in}
                  \addtolength{\hoffset}{#5}
                            \addtolength{\voffset}{#6}
                            \setlength{\oddsidemargin}{#1}
                            \setlength{\evensidemargin}{#2}
                            \setlength{\textwidth}{\paperwidth}
                  \addtolength{\textwidth}{-\oddsidemargin}
                  \addtolength{\textwidth}{-\evensidemargin}
                  \addtolength{\textwidth}{-\marginparsep}
                  \addtolength{\textwidth}{-\marginparwidth}
                            \setlength{\topmargin}{#3}
                            \setlength{\textheight}{\paperheight}
                  \addtolength{\textheight}{-\topmargin}
                  \addtolength{\textheight}{-\headheight}
                  \addtolength{\textheight}{-\headsep}
                  \addtolength{\textheight}{-\footskip}
                  \addtolength{\textheight}{-#4}}
\pageformat{2cm}{3cm}{25mm}{25mm}{1pt}{0pt}

\usepackage{ifthen}
\newboolean{article}
    \setboolean{article}{true}
\newboolean{report}
\newboolean{book}
\newboolean{letter}
\newboolean{german}
\newboolean{italian}
\newboolean{nobaselinestretch}
\newboolean{nosectionappendix}
\newboolean{oldtoc}
\newboolean{nosectionequn}
\newboolean{notheorem}

\ifthenelse{\boolean{german}}{
    \usepackage{german}}{}

\usepackage[latin1]{inputenc}

\usepackage{amsmath}
\usepackage{amssymb}
\usepackage[mathscr]{eucal}

\ifthenelse{\boolean{notheorem}}{}{
    \usepackage{theorem}}

\ifthenelse{\boolean{nobaselinestretch}}{}{
    \renewcommand{\baselinestretch}{1.25}}

\newenvironment{env}[2]{\begin{#1}#2\end{#1}}{}
    \newcommand{\beq}[1]{\begin{env}{equation}{#1}}
    \newcommand{\beqn}[1]{\begin{env}{equation*}{#1}}
    \newcommand{\bal}[1]{\begin{env}{align}{#1}}
    \newcommand{\baln}[1]{\begin{env}{align*}{#1}}
    \newcommand{\bga}[1]{\begin{env}{gather}{#1}}
    \newcommand{\bgan}[1]{\begin{env}{gather*}{#1}}
    \newcommand{\bflal}[1]{\begin{env}{flalign}{#1}}
    \newcommand{\bflaln}[1]{\begin{env}{flalign*}{#1}}
    \newcommand{\bmu}[1]{\begin{env}{multline}{#1}}
    \newcommand{\bmun}[1]{\begin{env}{multline*}{#1}}
    \newcommand{\bsp}[1]{\begin{env}{split}{#1}}

    \newcommand{\eeq}{\end{env}}
    \newcommand{\eeqn}{\end{env}}
    \newcommand{\eal}{\end{env}}
    \newcommand{\ealn}{\end{env}}
    \newcommand{\ega}{\end{env}}
    \newcommand{\egan}{\end{env}}
    \newcommand{\eflal}{\end{env}}
    \newcommand{\eflaln}{\end{env}}
    \newcommand{\emu}{\end{env}}
    \newcommand{\emun}{\end{env}}
    \newcommand{\esp}{\end{env}}

\newcommand{\lf}{\vspace{2ex}}

\renewcommand{\bf}[1]{\textbf{#1}}

\renewcommand{\sf}[1]{\textsf{#1}}

\renewcommand{\tt}[1]{\texttt{#1}}

\newcommand{\mbf}[1]{\mathbf{#1}}
\newcommand{\msf}[1]{\text{\small $\sf{#1}$}}

\newcommand{\cmc}[1]{\mathcal{#1}}
\newcommand{\eus}[1]{\mathscr{#1}}

\newcommand{\bb}[1]{\mathbb{#1}}

\newcommand{\nbd}[1]{$#1$\nobreakdash--}
\newcommand{\ol}[1]{\overline{#1}}

\newcommand{\ve}{\varepsilon}

\newcommand{\vp}{\varphi}

\newcommand{\abs}[1]{\left\lvert#1\right\rvert}
\newcommand{\norm}[1]{\left\lVert#1\right\rVert}

\newcommand{\Bnorm}[1]{\Bigl\lVert#1\Bigr\rVert}
\newcommand{\snorm}[1]{\norm{\smash{#1}}}
\newcommand{\sabs}[1]{\abs{\smash{#1}}}

\newcommand{\bfam}[1]{\bigl(#1\bigr)}
\newcommand{\Bfam}[1]{\Bigl(#1\Bigr)}
\newcommand{\AB}[1]{\langle#1\rangle}
\newcommand{\bAB}[1]{\bigl\langle#1\bigr\rangle}

\newcommand{\CB}[1]{\{#1\}}
\newcommand{\bCB}[1]{\bigl\{#1\bigr\}}
\newcommand{\BCB}[1]{\Bigl\{#1\Bigr\}}
\newcommand{\SB}[1]{[#1]}

\newcommand{\RO}[1]{[#1)}

\newcommand{\set}[2][]{
    \ifthenelse{\equal{#1}{}}{
        \CB{#2}}{
        \CB{#1~|~#2}}}
\newcommand{\bset}[2][]{
    \ifthenelse{\equal{#1}{}}{
        \bCB{#2}}{
        \bCB{#1~|~#2}}}
\newcommand{\Bset}[2][]{
    \ifthenelse{\equal{#1}{}}{
        \BCB{#2}}{
        \BCB{#1~\big|~#2}}}
\newcommand{\zero}{\CB{0}}

\DeclareMathOperator{\ls}{\normalfont\msf{span}}

\DeclareMathOperator{\cls}{\ol{\ls}}

\DeclareMathOperator{\id}{\normalfont\msf{id}}

\renewcommand{\ker}{\operatorname{\msf{ker}}}

\newcommand{\N}{\bb{N}}

\newcommand{\cA}{\cmc{A}}
\newcommand{\cB}{\cmc{B}}

\newcommand{\sB}{\eus{B}}

\newcommand{\sK}{\eus{K}}

\newcommand{\U}{\mbf{1}}

\ifthenelse{\boolean{nosectionequn}}{}{
    \numberwithin{equation}{section}
    }

\ifthenelse{\boolean{article}\or\boolean{letter}\or\boolean{nosectionequn}}{
    \setboolean{nosectionappendix}{true}}{}
\ifthenelse{\boolean{nosectionappendix}}{}{
    \renewcommand{\appendix}{
        \chapter*{\appendixname}
        \addcontentsline{toc}{chapter}{\appendixname}
        \renewcommand{\thesection}{\Alph{section}}
        \setcounter{section}{0}}}
   
\ifthenelse{\boolean{report}\or\boolean{book}}{
    }{}

\ifthenelse{\boolean{notheorem}}{}{
        \newcommand{\notename}{Note.}
        \newcommand{\mnname}{Mathematical note.}
        \newcommand{\enname}{End of the note.}
        \newcommand{\definame}{Definition.}
        \newcommand{\propname}{Proposition.}
        \newcommand{\lemname}{Lemma.}
        \newcommand{\exname}{Example.}
        \newcommand{\exername}{Exercise.}
        \newcommand{\remname}{Remark.}
        \newcommand{\obname}{Observation.}
        \newcommand{\thmname}{Theorem.}
        \newcommand{\corname}{Corollary.}
        \newcommand{\proofname}{Proof.}
    \ifthenelse{\boolean{german}}{
        \renewcommand{\mnname}{Mathematische Notiz.}
        \renewcommand{\enname}{Ende der Notiz.}
        \renewcommand{\exname}{Beispiel.}
        \renewcommand{\exername}{Übung.}
        \renewcommand{\remname}{Bemerkung.}
        \renewcommand{\obname}{Beobachtung.}
        \renewcommand{\thmname}{Satz.}
        \renewcommand{\corname}{Korollar.}
        \renewcommand{\proofname}{Beweis.}}{}
    \ifthenelse{\boolean{italian}}{
        \renewcommand{\mnname}{Nota matematica.}
        \renewcommand{\enname}{Fina della nota.}
        \renewcommand{\definame}{Definizione.}
        \renewcommand{\propname}{Proposizione.}
        \renewcommand{\exname}{Esempio.}
        \renewcommand{\exername}{Esercizio.}
        \renewcommand{\remname}{Nota.}
        \renewcommand{\obname}{Osservazione.}
        \renewcommand{\thmname}{Teorema.}
        \renewcommand{\corname}{Corollario.}
        \renewcommand{\proofname}{Dimostrazione.}

       \renewcommand{\appendixname}{Appendice}

       }{}
    \theoremheaderfont{\normalfont\bfseries}
    \theoremstyle{change}
        \theorembodyfont{\rmfamily}
            \newtheorem{emp}{}[section]
                \newcommand{\bemp}[1][]{
                    \begin{emp}\hskip-\labelsep\bf{#1}\hskip\labelsep}
                \newcommand{\eemp}{\end{emp}}
\newtheorem{itemp}[emp]{}
                \newcommand{\bitemp}[1][]{
                    \begin{itemp}\hskip-\labelsep\bf{#1}\hskip\labelsep\normalfont\itshape}
                \newcommand{\eitemp}{\end{itemp}}
            \newtheorem{note}[emp]{\notename}
                \newcommand{\bnote}{\begin{note}}
                \newcommand{\enote}{\end{note}}
            \newtheorem{mn}[emp]{\mnname}
                \newcommand{\bnm}{\begin{mn}~\begin{quotation}\renewcommand{\baselinestretch}{1}\small\noindent\ignorespaces}
                \newcommand{\enm}{\end{quotation}\hfill\bf{\enname}\end{mn}}
            \newtheorem{ex}[emp]{\exname}
                \newcommand{\bex}{\begin{ex}}
                \newcommand{\eex}{\end{ex}}
            \newtheorem{exer}[emp]{\exername}
                \newcommand{\bexer}{\begin{exer}}
                \newcommand{\eexer}{\end{exer}}
            \newtheorem{defi}[emp]{\definame}
                \newcommand{\bdefi}{\begin{defi}}
                \newcommand{\edefi}{\end{defi}}
            \newtheorem{rem}[emp]{\remname}
                \newcommand{\brem}{\begin{rem}}
                \newcommand{\erem}{\end{rem}}
            \newtheorem{ob}[emp]{\obname}
                \newcommand{\bob}{\begin{ob}}
                \newcommand{\eob}{\end{ob}}
        \theorembodyfont{\normalfont\itshape}
            \newtheorem{thm}[emp]{\thmname}
                \newcommand{\bthm}{\begin{thm}}
                \newcommand{\ethm}{\end{thm}}
            \newtheorem{prop}[emp]{\propname}
                \newcommand{\bprop}{\begin{prop}}
                \newcommand{\eprop}{\end{prop}}
            \newtheorem{cor}[emp]{\corname}
                \newcommand{\bcor}{\begin{cor}}
                \newcommand{\ecor}{\end{cor}}
            \newtheorem{lem}[emp]{\lemname}
                \newcommand{\blem}{\begin{lem}}
                \newcommand{\elem}{\end{lem}}
\newenvironment{empn}[1]{\lf\noindent\bf{#1}\ignorespaces\hskip\labelsep}{\lf}
		\newcommand{\bempn}[1]{\begin{empn}{#1}}
		\newcommand{\eempn}{\end{empn}}
		\newcommand{\bitempn}[1]{\begin{empn}{#1}\normalfont\itshape}
		\newcommand{\eitempn}{\end{empn}}
                \newcommand{\bnmn}{\begin{empn}{\mnname}~\begin{quotation}\renewcommand{\baselinestretch}{1}\small\noindent\ignorespaces}
                \newcommand{\enmn}{\end{quotation}\hfill\bf{\enname}\end{empn}}
		\newcommand{\bexn}{\begin{empn}{\exname}}
		\newcommand{\eexn}{\end{empn}}
		\newcommand{\bexern}{\begin{empn}{\exername}}
		\newcommand{\eexern}{\end{empn}}
		\newcommand{\bdefin}{\begin{empn}{\definame}}
		\newcommand{\edefin}{\end{empn}}
		\newcommand{\bremn}{\begin{empn}{\remname}}
		\newcommand{\eremn}{\end{empn}}
		\newcommand{\bobn}{\begin{empn}{\obname}}
		\newcommand{\eobn}{\end{empn}}

\newcommand{\qedsymbol}{~\rule[-0.35mm]{2mm}{2mm}}
    \newcounter{proof}[emp]
    \newenvironment{Proof}[1]{
        \vspace{1ex}
        \renewcommand{\item}[1][\stepcounter{proof}(\roman{proof})]%
            {##1\hskip\labelsep}
        \noindent\textsc{#1\hskip\labelsep}}{
        \nolinebreak\qedsymbol}
    \newcommand{\proof}[1][\proofname]{
        \begin{Proof}{#1}\ignorespaces}
    \newcommand{\qed}{\end{Proof}}
    \newcommand{\noqed}{
        \renewcommand{\qedsymbol}{}
        \end{Proof}}}
    \ifthenelse{\boolean{italian}}{
        \renewcommand{\proofname}{Dimostrazione.}}{}

\usepackage[varg]{txfonts}

\usepackage{stmaryrd}

\renewcommand{\thefootnote}{[\arabic{footnote}]}

\usepackage{hyperref}

\usepackage{color}

\def\blu{}

\begin{document}

\bibliographystyle{amsalpha}

\title{Kernels of Hilbert Module Maps:\\A Counterexample} 


\author{Jens Kaad and Michael Skeide{\renewcommand{\thefootnote}{}\footnote{MSC 2020: 46L08. Keywords: Hilbert modules, kernels, right linear bounded maps, orthogonal complement, hereditary subalgebra.}}
\setcounter{footnote}{0}
}

\date{January 2021}

\maketitle 

\vfill

\begin{abstract}
\noindent
Answering a long standing question, we give an example of a Hilbert module and a nonzero bounded right linear map having a kernel with trivial orthogonal complement. In particular, this kernel is different from its own double orthogonal complement.
%
%
\end{abstract}

\newpage

\section{Introduction}
The complement of a subset $S$ of a pre-Hilbert module $E$ is defined as
\beqn{
S^\perp
~:=~
\CB{x\in E\colon\AB{S,x}=\zero}.
}\eeqn
 It is well known that, given a subset $S$ of a Hilbert space $H$, the double complement $S^{\perp\perp}=(S^\perp)^\perp$ can be computed as
 \beqn{
 S^{\perp\perp}
 ~=~
 \cls S.
 }\eeqn
Since a bounded linear functional that vanishes on $S$ also vanishes on $\cls S$, it vanishes on $S^{\perp\perp}$. In particular, the kernel of a bounded linear functional on a Hilbert space coincides with its double complement. (In some papers this is phrased, saying the kernel is orthogonally complete.)

The same sort of reasoning shows that the kernel of a bounded \nbd{\cB}linear map from a von Neumann \nbd{\cB}module (or \nbd{W^*}module%
\footnote{
A \nbd{W^*}module is a self-dual Hilbert module over a \nbd{W^*}algebra; this definition is probably due to Baillet, Denizeau, and Havet \cite{BDH88}. A von Neumann module is a Hilbert module over a von Neumann algebra that is strongly closed in a uniquely associated operator space; see Skeide \cite{Ske00b}.
}%
) $E$ into $\cB$ coincides with its double complement. 

But, is it true for Hilbert modules? The question in its purest form is:
\begin{itemize}
\item[]
Given a bounded right linear map $\Phi$ from a Hilbert \nbd{\cB}module $E$ to $\cB$, does $\Phi(S)=\zero$ for a subset $S\subset E$ with $S^\perp=\zero$ imply that $\Phi=0$?
\end{itemize}
In its most general form, the question is:
\begin{itemize}
\item[]
Given a bounded right linear map $a$ from a Hilbert \nbd{\cB}module $E$ to a Hilbert \nbd{\cB}module $E'$, do we always have $\ker a=(\ker a)^{\perp\perp}$?
\end{itemize}
These two questions have the same answer. Moreover, for special cases the answer is yes, for instance, if the map $a$ is required to be adjointable, or if {$\,\,\cls S\cB\,$} is an ideal submodule (Bakic and Guljas \cite{BaGu02}) or a closed ternary ideal (Skeide \cite{Ske18p1}). 

\lf
Several publications point out pleasant consequences of a positive answer to the above questions. (See Frank \cite{FraM02}, Bhat and Skeide \cite[Footnotes 1-3]{BhSk15}, and \cite[Footnote 7]{Ske18p1}.) But, is it true? The scope of this short note is to illustrate by an explicit counterexample that, unfortunately, the general answer is no:

\bthm \label{mthm}
There exist a Hilbert \nbd{\cB}module $E$, a bounded \nbd{\cB}linear map $\Phi\colon E\rightarrow\cB$, and a closed submodule $F$ such that
\beqn{
E
~\ne~
\ker\Phi
~\supset~
F
\text{~~~~~~and~~~~~~}
F^\perp
~=~
\zero.
}\eeqn
In particular, $\ker\Phi^{\perp\perp}\supset F^{\perp\perp}=E\ne\ker\Phi$.
\ethm

And the answer is no under the simplest imaginable circumstances: $E$ is a standard Hilbert module over a separable unital commutative \nbd{C^*}algebra and $F$ is a full submodule.%
\footnote{
This does not mean that the statements in \cite{FraM02,BhSk15,Ske18p1} cannot be true. In fact, Shalit and Skeide \cite[Appendix B]{ShaSk10p} prove, by a different method, a result that {\blu would have followed from a positive answer to the above questions}; this is what made the second named author start to {\blu think about them}.}

\lf
Before diving into the example in the next section, let us briefly discuss why, we think, the statement we disprove in this note is suspicious right from the beginning.

Hilbert modules are in particular Banach spaces and, of course, share all the good properties of general Banach spaces. But as far as their geometric properties are concerned, frequently Hilbert modules behave much more like pre-Hilbert spaces than like Hilbert spaces. (Closed submodules need not be complemented; bounded right linear maps need not have an adjoint; in fact, an isometry has an adjoint if and only if its range is complemented. All these failures, in the end, go back to the fact that a Hilbert module need not be self-dual; and our example is no exception.) And for pre-Hilbert spaces the answer is no. We give a quick example, which we owe to Orr Shalit. This example both indicates how to find a counterexample and explains what goes wrong in existing attempts to find a positive answer to the above questions.

\bex
Let $H$ be a pre-Hilbert space with orthonormal Hamel basis $\bfam{e_n}_{n\in\N}$ and let $S:=\CB{e_n-2e_{n+1}}$. Then $S^\perp$ is $\zero$ in $H$. But, $S^\perp$ is not $\zero$ in the completion $\ol{H}$. Therefore, for every element $\vp\ne0$ in the complement in $\ol{H}$, the functional $\AB{\vp,\bullet}$ vanishes on $S$, but not on $H$.
\eex

Under the passage from $H$ to $\ol{H}$, the set $S$ loses the property to be \nbd{0}complemented. So, knowing the answer is yes for Hilbert spaces, does not help proving the same for pre-Hilbert spaces by embedding the latter into the former. Based on the observation that for von Neumann modules the answer is affirmative, there have been attempts to prove the same for Hilbert modules by embedding the latter into the former. Specifically for the embedding into the bidual it follows from a careful analysis of results by Akemann and Pedersen \cite{AkPe73} (see \cite[Section 3.11]{Ped79} in Pedersen's book) that this can never work.

We end this note by providing an application to general \nbd{C^*}algebra theory concerning representations of hereditary \nbd{C^*}subalgebras.
%


\section{The counterexample}
Let $\cB$ be a \nbd{C^*}algebra. Recall that for a set $S$ the standard Hilbert \nbd{\cB}module over $S$ is defined as
\beqn{
\cB^S
~=~
\BCB{\bfam{b_s}_{s\in S}\colon\sum_{s\in S}b_s^*b_s\text{~exists as norm limit over the finite subsets of $S$}}.
}\eeqn
If $\cB$ is unital, then $\cB^S$ has an orthonormal basis $\bfam{e_s}_{s\in S}$ defined by $e_s:=\bfam{\delta_{s,s'}\U}_{s'\in S}$. It is well known that the formula
\beqn{
\beta_s^*
~=~
\Phi(e_s)
}\eeqn
establishes (for unital $\cB$) a one-to-one correspondence between bounded right linear maps $\Phi\colon\cB^S\rightarrow\cB$ and families $\bfam{\beta_s}_{s\in S}$ ($\beta_s\in\cB$) such that {\blu there exists a constant $M\ge 0$} with
\beqn{
\Bnorm{\sum_{s\in S'}\beta_s^*\beta_s}
~\leq~
M
}\eeqn
for all finite subsets $S'\subset S$. Indeed, {\blu suppose that we are given a family $\bfam{\beta_s}_{s\in S}$ for which such a constant $M \ge 0$ exists, then} putting $B_{S'}:=\bfam{\beta_s\chi_{S'}(s)}_{s\in S}\in\cB^S$, from $\snorm{B_{S'}B_{S'}^*}=\norm{\AB{B_{S'},B_{S'}}}$ we see that the net {\blu $\bfam{B_{S'}B_{S'}^*}_{S'}$} of rank-one operators is bounded by $M$. Clearly, we have $\lim_{S'}\AB{B_{S'},e_s}=\Phi(e_s)$. This means, the bounded net {\blu $\bfam{\AB{B_{S'},\bullet}}_{S'}$} converges strongly to $\Phi$ on the dense subset of finite \nbd{\cB}linear combinations of the basis vectors $e_s$, hence, everywhere. See Manuilov and Troitsky \cite[Section 2.5]{MaTr05} for a different proof for sequences. By \bf{the} standard Hilbert \nbd{\cB}module, we mean $\cB^\N$, frequently also written as $H_\cB$ or $\cB^\infty$.

For our example we fix $\cB:=C\SB{0,1}$. For $q\in (0,1]$ choose a sequence of elements $(\psi_{q,m})_{m \in \N}$ in $\cB$ such that
\beqn{
\sum_{m=1}^\infty\sabs{\psi_{q,m}(t)}^2
~=~
\chi_{\RO{0,q}}(t)
}\eeqn
pointwise for $t\in\SB{0,1}$. (For instance, with $f_{q,m}\colon t\mapsto\max(0,\min(1,m(q-t)))$ and $f_{q,0}=0$, putting $\psi_{q,m}:=\sqrt{f_{q,m}-f_{q,m-1}}$ will do.) It follows that $\Psi_q=\bfam{\psi_{q,m}}_{m\in\N}$ defines a bounded right linear map $\cB^\N\rightarrow\cB$. (Since $\bfam{\psi_{q,m}}_{m\in\N}$ is not an element of $\cB^\N$, the map $\Psi_q$ is not adjointable.)

Next, we put
\beqn{
E
~:=~
\cB\oplus(\cB^\N)^\infty,
}\eeqn
and, with the set $A:=\zero \cup ( \N\times\N )$, we shall write $E$ as
\beqn{
E
~=~
\cB^A.
}\eeqn
The index $0\in A$ refers to the separate copy of $\cB$. The second index $m$ of $(n,m)\in A$ refers to the $m$th element of the sequence in $\cB^\N$, while the first index $n$ refers to the number of the copy of $\cB^\N$ we are talking about.

Let $\bfam{q_n}_{n\in\N}$ be a sequence of elements in $(0,1]$. By the assignment
\beqn{
\Phi
\colon
\bfam{b_a}_{a\in A}
~\longmapsto~
b_0+\sum_{n=1}^\infty\textstyle {2^{-n}} \Psi_{q_n}\Bfam{\bfam{b_{(n,m)}}_{m\in\N}},
}\eeqn
we define a nonzero right linear bounded map. In fact, $\Phi \colon E \rightarrow \cB$ is surjective.

If we define the elements
\beqn{
\zeta_{k,\ell}
~:=~
\bfam{\delta_{a,0}\textstyle {2^{-k} } \psi_{q_k,\ell}^*-\delta_{a,(k,\ell)}\U}_{a\in A}
}\eeqn
in $E$, then $\Phi(\zeta_{k,\ell})=0$. So, for $S:= \CB{\zeta_{k,l}\colon k,\ell\in\N}$ the closed submodule $F:=\cls S\cB$ is contained in $\ker\Phi$.

For the punchline, we are done proving Theorem \ref{mthm}, if we show that the sequence $(q_n)_{n \in \N}$ can be chosen such that $S^\perp=\zero$.

\blem
Suppose that the set $\CB{q_n\colon n\in\N}$ is dense $\SB{0,1}$. Then $S^\perp=\zero$.
\elem

\proof
Let $\bfam{b_a}_{a\in A} {\in \cB^A}$ and suppose that $\bAB{\zeta_{n,m},\bfam{b_a}_{a\in A}}=0$ for all $n,m$. Then
\beqn{ \tag{$*$}
{2^{-n}}\psi_{q_n,m}b_0-{b_{n,m}}
~=~
0.
}\eeqn
That is, each $b_{n,m}$ is determined by $b_0$. In particular, if $b_0=0$, then so is $b$.

{ Suppose that $b_0\ne0$ and put $b_{n,m}$ as determined by $(*)$. Then by density of the {\blu subset $\CB{q_n\colon n\in\N} \subset \SB{0,1}$} and continuity of $b_0$ there are $n\in\N$, $\ve>0$, and $d>0$ such that $\abs{b_0(t)}\ge d$ for all $t\in(q_n-\ve,q_n+\ve)\cap\SB{0,1}$. We know that
{\beqn{
f(t)
~:=~
\sum_{m = 1}^\infty \sabs{b_{n,m}(t)}^2 
~=~
\sum_{m = 1}^\infty 4^{-n} \sabs{\psi_{q_n,m}(t) b_0(t)}^2
~=~ 4^{-n} \sabs{b_0(t)}^2 \chi_{[0,q_n)}(t)
}\eeqn 
for all $t \in [0,1]$. Since $\chi_{[0,q_n)}(t)$ has a discontinuity at the point $q_n \in (0,1]$ and since $\abs{b_0(t)}\ge d$ on a neighbuorhood of $q_n$, also $f$ is not continuous at $q_n$. It follows that} $\sum_{m=1}^\infty\sabs{b_{n,m}}^2$ does not converge uniformly in {$\SB{0,1}$, so $\bfam{b_a}_{a\in A}\notin \cB^A$.} Consequently, if $\bfam{b_a}_{a\in A}\in S^\perp\subset\cB^A$, then $b_0=0$ hence $\bfam{b_a}_{a\in A}=0$.}\qed

\section{An application}
The following result about representations of hereditary \nbd{C^*}subalgebras might be of interest for general \nbd{C^*}algebra theory.

\bthm
There is a \nbd{C^*}algebra $\cA$ with hereditary \nbd{C^*}subalgebra $\cB$ separating the points of $\cA$ (that is, $a\cB=\zero$ implies $a=0$) such that there is no faithful representation $\pi$ of $\cA$ such that $\pi(\cB)$ acts nondegenerately on the representation space.
\ethm

\proof
In the example from the previous section, the \nbd{C^*}algebra $\sK(E)$ and its hereditary subalgebra $\sK(F)$ do the job. Indeed, the Hilbert \nbd{\sK(E)}module $\sK(E)$ may be identified with the (internal) tensor product $E\odot E^*$, and it contains $\sK(E,F)=F\odot E^*$ as zero-complemented submodule. $\Phi$ gives rise to a nonzero bounded right-linear map $\Phi\odot\id_{E^*}\colon E\odot E^*\rightarrow E^*$ that vanishes on $F\odot E^*$.

Let $\pi\colon\sK(E)\rightarrow\sB(H)$ be a representation. Then $E\odot E^*\odot H=\cls\pi(\sK(E))H$ via the identification $x'\odot x^*\odot h\mapsto\pi(x'x^*)h$. So, suppose $\pi$ is nondegenerate, so that $E\odot E^*\odot H=H$. The map $\Phi$ induces the bounded map $\Phi\odot{\blu \id_{E^*}}\odot\id_H$ from $H=E\odot E^*\odot H$ to $E^*\odot H$. Suppose $\pi$ is faithful, so that $\Phi\odot{\blu \id_{E^*}}\odot\id_H\ne0$ (because, otherwise, $\Phi=0$). Since $\Phi\odot{\blu \id_{E^*}}\odot\id_H$ vanishes on $F\odot E^*\odot H=\cls\pi(F\odot E^*)H {\blu \supset} \cls\pi(F\odot F^*)H$, the representation $\pi$ restricted to $\sK(F)=F\odot F^*$ cannot be nondegenerate (because, otherwise, $\Phi$ would be $0$).\qed

\newpage

 \noindent
 \bf{Acknowledgments:} This little note is an output of a talk by the second author at the ``Hilbert \nbd{C^*}Modules Online Weekend'' from December 5/6, 2020, presenting the open problem {which was} then answered by the first author. We are deeply grateful to the organizers of the Workshop Michael Frank, Vladimir Manuilov, and Evgenij Troitsky for having made this possible. We also acknowledge useful discussions with Michael Frank, Boris Guljas, and Orr Shalit.

\setlength{\baselineskip}{2.5ex}


\newcommand{\Swap}[2]{#2#1}\newcommand{\Sort}[1]{}
\providecommand{\bysame}{\leavevmode\hbox to3em{\hrulefill}\thinspace}
\providecommand{\MR}{\relax\ifhmode\unskip\space\fi MR }
\providecommand{\MRhref}[2]{%
  \href{http://www.ams.org/mathscinet-getitem?mr=#1}{#2}
}
\providecommand{\href}[2]{#2}

\lf\noindent
Jens Kaad:
{\small\itshape Department of Mathematics and Computer Science, The University of Southern Denmark, Campusvej 55, DK-5230 Odense M, Denmark, E-mail: \href{mailto:kaad@imada.sdu.dk}{\tt{kaad@imada.sdu.dk}}}\\
{\small{\itshape Homepage: \url{https://portal.findresearcher.sdu.dk/en/persons/kaad}}}

\lf\noindent
Michael Skeide:
{\small\itshape Dipartimento di Economia, Università degli Studi del Molise, Via de Sanctis, 86100 Campobasso, Italy, E-mail: \href{mailto:skeide@unimol.it}{\tt{skeide@unimol.it}}}\\
{\small{\itshape Homepage: \url{http://web.unimol.it/skeide/}}}


\end{document}